\documentclass[12pt,oneside,reqno]{amsart}
\usepackage{txfonts}
\usepackage{bbm}
\usepackage{amsmath}
\usepackage{graphicx}
\usepackage{mathrsfs}
\usepackage{stmaryrd}
\usepackage{amsfonts}
\usepackage{enumerate,amsmath,amssymb,amsthm}

\numberwithin{equation}{section}

\newcommand{\be}{\begin{eqnarray}}
\newcommand{\ee}{\end{eqnarray}}
\newcommand{\ce}{\begin{eqnarray*}}
\newcommand{\de}{\end{eqnarray*}}
\newtheorem{theorem}{Theorem}[section]
\newtheorem{lemma}[theorem]{Lemma}
\newtheorem{remark}[theorem]{Remark}
\newtheorem{definition}[theorem]{Definition}
\newtheorem{proposition}[theorem]{Proposition}
\newtheorem{Examples}[theorem]{Example}
\newtheorem{corollary}[theorem]{Corollary}

\def\eps{\varepsilon}

\def\a{\alpha}

\def\p{\partial}

\def\[{{\Big[}}
\def\]{{\Big]}}
\def\<{{\langle}}
\def\>{{\rangle}}
\def\({{\Big(}}
\def\){{\Big)}}

\def\bx{{\mathbf{x}}}

\def\sgn{\mbox{\rm sgn}}

\def\dif{{\mathord{{\rm d}}}}

\def\min{{\mathord{{\rm min}}}}

\def\no{\nonumber}
\def\={&\!\!=\!\!&}
\def\bt{\begin{theorem}}
\def\et{\end{theorem}}
\def\bl{\begin{lemma}}
\def\el{\end{lemma}}
\def\br{\begin{remark}}
\def\er{\end{remark}}

\def\bd{\begin{definition}}
\def\ed{\end{definition}}
\def\bp{\begin{proposition}}
\def\ep{\end{proposition}}
\def\bc{\begin{corollary}}
\def\ec{\end{corollary}}
\def\bx{\begin{Examples}}
\def\ex{\end{Examples}}

\def\cB{{\mathcal B}}

\def\cE{{\mathcal E}}

\def\cJ{{\mathcal J}}

\def\cM{{\mathcal M}}

\def\cP{{\mathcal P}}

\def\cT{{\mathcal T}}

\def\mD{{\mathbb D}}
\def\mE{{\mathbb E}}

\def\mI{{\mathbb I}}

\def\mM{{\mathbb M}}
\def\mN{{\mathbb N}}

\def\mR{{\mathbb R}}

\def\mU{{\mathbb U}}

\def\sF{{\mathscr F}}

\def\sL{{\mathscr L}}
\def\sM{{\mathscr M}}

\def\sU{{\mathscr U}}

\def\geq{\geqslant}
\def\leq{\leqslant}

\def\div{\mathord{{\rm div}}}

\allowdisplaybreaks

\begin{document}

\title{Degenerate Irregular SDEs with Jumps and Application to Integro-Differential Equations of Fokker-Planck type}

\date{}
\author{Xicheng Zhang}


\dedicatory{
School of Mathematics and Statistics,
Wuhan University, Wuhan, Hubei 430072, P.R.China,\\
Email: XichengZhang@gmail.com
 }

\begin{abstract}
We investigate stochastic differential equations with jumps and irregular coefficients, and obtain the existence and uniqueness of
generalized stochastic flows.
Moreover, we also prove the existence and uniqueness of $L^p$-solutions or measure-valued solutions
for second order integro-differential equation of Fokker-Planck type.

\end{abstract}

\maketitle
\rm

\section{Introduction}

Recently, there are increasing interests to extend the classical DiPerna-Lions theory \cite{Di-Li} about ordinary differential equations (ODE)
with Sobolev coefficients to the case of stochastic differential equations (SDE)
(cf. \cite{Le-Li,Le-Li1,Fi, Zh1,Zh3, Zh2,Fa-Lu-Th,Li-Lu}).
In \cite{Fi}, Figalli first extended the DiPerna-Lions theory to
SDE in the sense of martingale solutions by using analytic tools and solving deterministic Fokker-Planck equations.
In \cite{Le-Li}, Le Bris and Lions studied the almost everywhere stochastic flow of SDEs
with constant diffusion coefficients, and in \cite{Le-Li1}, they also gave an outline for proving
the pathwise uniqueness for SDEs with irregular coefficients by studying the corresponding Fokker-Planck equations with
irregular coefficients. In \cite{Zh1} and \cite{Zh2}, we extended DiPerna-Lions' result to the case of SDEs  by using Crippa and
De Lellis' argument \cite{Cr-De-Le}, and obtained the existence and uniqueness of generalized stochastic flows for SDEs with irregular coefficients
(see also \cite{Fa-Lu-Th} for some related works).
Later on, Li and Luo \cite{Li-Lu} extended Ambrosio's result \cite{Am} to the case of SDEs with BV drifts and smooth diffusion
coefficients by transforming the SDE to an ODE. Moreover, a limit theorem for SDEs with discontinuous coefficients approximated by ODEs
was also obtained in \cite{Re-Zh}.

In this paper we are concerned with the following SDEs in $[0,1]\times\mR^d$ with jumps:
\begin{align}
\dif X_t=b_t(X_t)\dif t+\sigma_t(X_t)\dif W_t
+\int_{\mR^d\setminus\{0\}}f_t(X_{t-},y)\tilde N(\dif y,\dif t),\label{SDE01}
\end{align}
where $b:[0,1]\times\mR^d\to\mR^d$, $\sigma:[0,1]\times\mR^d\to\mR^d\times\mR^d$ and
$f:[0,1]\times\mR^d\times\mR^d\to\mR^d$ are measurable functions, $(W_t)_{t\in[0,1]}$
is a $d$-dimensional Brownian motion
and $N(\dif y,\dif t)$ is a Poisson random measure in $\mR^d\setminus\{0\}$ with intensity measure $\nu_t(\dif y)\dif t$,
$\tilde N(\dif y,\dif t):=N(\dif y,\dif t)-\nu_t(\dif y)\dif t$ is the compensated Poisson random measure.
The aim of the present paper is to extend the results in \cite{Zh1} to the above jump SDEs with Sobolev drift $b$
and Lipschitz $\sigma,f$.

Let us now describe the motivation. Suppose that $f_t(x,y)=y$.
Let $\sL$ be the generator of SDE (\ref{SDE01}) (a second order integro-differential operator) given as follows:
for $\varphi\in C^\infty_b(\mR^d)$, smooth function with bounded derivatives of all orders,
$$
\sL_t\varphi(x):=\frac{1}{2} a^{ij}_t(x)\p_i\p_j \varphi(x)+b^i_t(x)\p_i \varphi(x)+
\int_{\mR^d\setminus\{0\}}[\varphi(x+y)-\varphi(x)-y^i\p_i \varphi(x)]\nu_t(\dif y),
$$
where $a^{ij}_t(x):=\sum_k\sigma^{ik}_t(x)\sigma^{jk}_t(x)$, and we have used that the repeated
indices in a product is summed automatically, and this convention will be in forced throughout the present paper.
Here, we assume that for any $p\geq 1$,
\begin{align}
\int^1_0\!\!\!\int_{\mR^d\setminus\{0\}}|y|^2(1+|y|^2)^p\nu_s(\dif y)\dif s<+\infty.\label{BB3}
\end{align}
Let $X_t$ be a solution of SDE (\ref{SDE01}).
The law of $X_t$ in $\mR^d$ is denoted by $\mu_t$. Then by
It\^o's formula (cf. \cite{Ik-Wa} or \cite{Ap}), one sees that $\mu_t$ solves the following second order partial integro-differential equation (PIDE) 
of Fokker-Planck type in the distributional sense:
\begin{align}
\p_t\mu_t=\sL^*_t\mu_t,\label{FP}
\end{align}
subject to the initial condition:
\begin{align}
\lim_{t\downarrow 0}\mu_t=\mbox{Law of $X_0$ in the sense of weak convergence,}\label{Ini}
\end{align}
where $\sL^*_t$ is the adjoint operator of $\sL_t$ formally given by
$$
\sL^*_t\mu:=\frac{1}{2} \p_i\p_j (a^{ij}_t(x)\mu)-\p_i(b^i_t(x)\mu)+
\int_{\mR^d\setminus\{0\}}[\tau_{y}\mu-\mu+y^i\p_i\mu]\nu_t(\dif y),
$$
where for a probability measure $\mu$ in $\mR^d\setminus\{0\}$ and $y\in\mR^d$, $\tau_y\mu:=\mu(\cdot-y)$.
More precisely, for any $\varphi\in C^\infty_b(\mR^d)$,
\begin{align}
\p_t\<\mu_t,\varphi\>=\<\mu_t,\sL_t\varphi\>,\label{Dis}
\end{align}
where $\<\mu_t,\varphi\>:=\int_{\mR^d}\varphi(x)\mu_t(\dif x)$. If $b$ and $\sigma$ are not continuous,
in order to make sense for (\ref{Dis}), one needs to  at least assume that
$$
\int^1_0\!\!\!\int_{\mR^d}(|b_t(x)|+|a_t(x)|)\mu_t(\dif x)\dif t<+\infty.
$$
The following two questions are our main motivations of this paper:

($1^o$) Under what less conditions on the coefficients and in what spaces or senses does the uniqueness
for PIDE (\ref{FP})-(\ref{Ini}) hold?

($2^o$) If the initial distribution $\mu_0$ has a density with respect to the Lebesgue measure, does $\mu_t$ have a density
 with respect to the Lebesgue measure for any $t\in(0,1]$?

When there is no jump part and the diffusion coefficient is non-degenerate, in \cite{Bo-Da-Ro-St} the authors have
already given rather weak conditions for the uniqueness of measure-valued solutions based upon the Dirichlet form theory. In \cite{Fi}, Figalli
also gave some other conditions for the uniqueness of $L^1\cap L^\infty$-solutions by proving a maximal principle.
In \cite{Ro-Zh}, using a representation formula for the solutions of PDE (\ref{FP}) proved in \cite{Fi},
which is originally proved by Ambrosio \cite{Am} for continuity equation,
we gave different conditions for the uniqueness of measure-valued solutions and $L^p$-solutions to
second order degenerated Fokker-Planck equations. However, to the best of the author's knowledge, there are few results
on the integro-differential equation of Fokker-Planck type. The {\it non-local}
character of the operator $\sL$ causes some new difficulties to analyze by the classical tools.

For answering the above two questions to equation (\ref{FP}), we shall use a purely probabilistic approach.
The first step is to extend the almost everywhere stochastic flow in \cite{Le-Li,Zh1, Zh2} to SDE (\ref{SDE01})
so that we can solve the above question ($2^o$). In this extension, we need to carefully treat the jump size.
Since even in the linear case, if one does not
make any restriction on the jump, the law of the solution would not be absolutely continuous
with respect to the Lebesgue measure (cf. \cite[p.328, Example]{Pr}). The next step is to prove a
representation formula for the solution of (\ref{FP}) as in \cite[Theorem 2.6]{Fi}. This will lead to the uniqueness
of PIDE (\ref{FP}) by proving the pathwise uniqueness of SDE (\ref{SDE01}).

This paper is organized as follows: In Section 2, we collect some well known facts for later use. In Section 3,
we study the smooth SDEs with jumps, and prove an a priori estimate about the Jacobi determinant of $x\mapsto X_t(x)$.
In Section 4, we prove the existence and uniqueness of
almost everywhere or generalized stochastic flows for SDEs with jumps and rough drifts. In Section 5, the application to second order
integro-differential equations of Fokker-Planck type is presented.

\section{Preliminaries}

Throughout this paper we assume that $d\geq 2$.
Let $\mM_{d\times d}$ be the set of all $d\times d$-matrices. We need the following simple lemma about the
differentials of determinant function.
\bl\label{Le1}
Let $A=(a_{ij}), B=(b_{ij})\in\mM_{d\times d}$. Then the first and second order derivatives
of the determinant function $\det:\mM_{d\times d}\to\mR$ are given by
\begin{align}
(\nabla\det)(A)(BA):=\frac{\dif }{\dif t}\det(A+tBA)|_{t=0}=\det(A)\mathrm{tr}(B)\label{Ep3}
\end{align}
and
\begin{align}
(\nabla^2\det)(A)(BA, BA):=\frac{\p^2 }{\p t\p s}\det(A+tBA+sBA)|_{s=t=0}
=\det(A)\sum_{i,j}[b_{ii}b_{jj}-b_{ij}b_{ji}].\label{Ep4}
\end{align}
Moreover, if $|b_{ij}|\leq\a$ for all $i,j$, then
\begin{align}
|\det(\mI+B)-1-\mathrm{tr}(B)|\leq d! d^2\a^2(1+\a)^{d-2}.\label{Ep5}
\end{align}
\el
\begin{proof}
Notice that
$$
\det(A+tBA)=\det(A)\det(\mI+tB)
$$
and
$$
\det(A+tBA+sBA)=\det(A)\det(\mI+(t+s)B).
$$
Formulas (\ref{Ep3}) and (\ref{Ep4}) are easily derived from the definition
\begin{align}
\det(\mI+tB):=\sum_{\sigma\in S_d}\sgn(\sigma)\prod_{i=1}^d(1_{i\sigma(i)}+tb_{i\sigma(i)}),\label{Ep6}
\end{align}
where $S_d$ is the set of all permutations of $\{1,2,\cdots, d\}$ and $\sgn(\sigma)$ is the sign of $\sigma$.

As for (\ref{Ep5}), let $h(t):=\det(\mI+tB)$, then $h'(0)=\mathrm{tr}(B)$ and
$$
\det(\mI+B)-1-\mathrm{tr}(B)=\int^1_0\!\!\!\int^t_0h''(s)\dif s\dif t=\int^1_0(1-s)h''(s)\dif s.
$$
Estimate (\ref{Ep5}) now follows from (\ref{Ep6}).
\end{proof}
The following result is taken from \cite[Theorem 6]{Pr-Sh}.
\bt\label{Th1}  Let $M$ be a locally square integrable martingale such that $\Delta M>-1$ a.s. Let
$\cE(M)$ be the Dol\'eans-Dade exponential defined by
$$
\cE(M)_{t}:=\exp\Big\{M_t-\frac{1}{2}\<M^{\mathrm{c}}\>_t\Big\}\times
\prod_{0<s\leq t}(1+\Delta M_{s})e^{-\Delta M_{s}}.
$$
If for some $T>0$,
$$
\mE\Big[\exp\Big\{\frac{1}{2}\<M^{\mathrm{c}}\>_T+\<M^{\mathrm{d}}\>_T\Big\}\Big]<\infty,
$$
where $M^{\mathrm{c}}$ and $M^{\mathrm{d}}$ are respectively continuous and purely discontinuous martingale parts of $M$, then $\cE(M)$
is a martingale on $[0,T]$.
\et

In Sections 3 and 4, we shall deal with the general Poisson point process. Below we introduce
some necessary spaces and processes. Let $(\Omega,\sF,P;(\sF_t)_{t\geq 0})$ be a complete filtered
probability space and $(\mU,\sU)$ a measurable space. Let $(W(t))_{t\geq 0}$ be a $d$-dimensional standard
($\sF_t$)-adapted Brownian motion and $(p_t)_{t\geq 0}$ an ($\sF_t$)-adapted Poisson point process with values in $\mU$ and with
intensity measure $\nu_t(\dif u)\dif t$, a $\sigma$-finite measure on $[0,1]\times\mU$ (cf. \cite{Ik-Wa}).
Let $N(\dif u,(0,t])$ be the counting measure of $p_t$, i.e., for any $\Gamma\in\sU$,
$$
N(\Gamma,(0,t]):=\sum_{0<s\leq t}1_\Gamma(p_s).
$$
The compensated Poisson random measure of $N$ is given by
$$
\tilde{N}(\dif u,(0,t]):=N(\dif u,(0,t])-\int^t_0\nu_s(\dif u)\dif s.
$$
We remark that for $\Gamma\in\sU$ with $\int^t_0\nu_s(\Gamma)\dif s<+\infty$, the random variable $N((0,t],\Gamma)$
obeys the Poisson distribution with parameter $\int^t_0\nu_s(\Gamma)\dif s$.

Below, the letter $C$ with or without subscripts
will denote a positive constant whose value is not important and may change in different occasions.
Moreover, all the derivatives, gradients and divergences are taken in the distributional sense.

The following lemma is a generalization of \cite[Proposition 1.12, p. 476]{Re-Yo}
(cf. \cite[Lemma A.2]{Qi-Zh}).
\bl\label{Le3}
Let $L:\mU\to\mR$ be a measurable function satisfying that $|L(u)|\leq C$ and $\int^1_0\!\!\!\int_\mU L(u)^2\nu_s(\dif u)\dif s<+\infty$.
Then for any $t>0$,
$$
\mE\exp\left\{\sum_{0<s\leq t}L(p_s)^2\right\}
=\exp\left\{\int^t_0\!\!\!\int_{\mU}(e^{L(u)^2}-1)\nu_s(\dif u)\dif s\right\}<+\infty.
$$
\el

We also need the following technical lemma (cf. \cite[Lemma 3.4]{Zh2}).
\bl\label{Le4}
Let $\mu$ be a locally finite measure on $\mR^d$ and
$(X_n)_{n\in\mN}$ be a family of random fields on $\Omega\times\mR^d$.
Suppose that $X_n$ converges to $X$ for $P\otimes\mu$-almost all $(\omega,x)$, and for some $p\geq 1$,
there is a constant $K_p>0$ such that for any nonnegative
measurable function $\varphi\in L^p_\mu(\mR^d)$,
\begin{align}
\sup_n\mE\int_{\mR^d}\varphi(X_n(x))\mu(\dif x)\leq K_p\|\varphi\|_{L^p_\mu}.\label{Lp1}
\end{align}
Then we have:

(i). For any nonnegative measurable function $\varphi\in L^p_\mu(\mR^d)$,
\begin{align}
\mE\int_{\mR^d}\varphi(X(x))\mu(\dif x)\leq K_p\|\varphi\|_{L^p_\mu}.\label{Ps2}
\end{align}

(ii). If $\varphi_n$ converges to $\varphi$ in $L^p_\mu(\mR^d)$, then for any $N>0$,
\begin{align}
\lim_{n\to\infty}\mE\int_{|x|\leq N}|\varphi_n(X_n(x))-\varphi(X(x))|\mu(\dif x)=0.\label{Ps22}
\end{align}
\el

Let $\varphi$ be a locally integrable function
on $\mR^d$. For every $R>0$, the local maximal function is defined by
$$
\cM_R \varphi(x):=\sup_{0<r<R}\frac{1}{|B_r|}\int_{B_r}\varphi(x+y)\dif y
=:\sup_{0<r<R}\fint_{B_r} \varphi(x+y)\dif y,
$$
where $B_r:=\{x\in\mR^d: |x|<r\}$ and $|B_r|$ denotes the volume of $B_r$.
The following result can be found in \cite[p.143, Theorem 3]{Ev-Ga} and \cite[Appendix A]{Cr-De-Le}.
\bl
(i) (Morrey's inequality)
Let $\varphi\in L^1_{loc}(\mR^d)$ be such that $\nabla\varphi\in L^q_{loc}(\mR^d)$ for some $q>d$. Then
there exist $C_{q,d}>0$ and a negligible set $A$ such that
for all $x,y\in A^c$ with $|x-y|\leq R$,
\begin{align}
|\varphi(x)-\varphi(y)|&\leq C_{q,d}\cdot |x-y|\cdot\left(\fint_{B_{|x-y|}}
|\nabla \varphi|^q(x+z)\dif z\right)^{1/q}\no\\
&\leq C_{q,d}\cdot |x-y|\cdot (\cM_R|\nabla \varphi|^q(x))^{1/q}.\label{Es02}
\end{align}
(ii) Let $\varphi\in L^1_{loc}(\mR^d)$ be such that $\nabla\varphi\in L^1_{loc}(\mR^d)$.
Then there exist $C_d>0$ and a negligible set $A$ such that
for all $x,y\in A^c$ with $|x-y|\leq R$,
\begin{align}
|\varphi(x)-\varphi(y)|\leq C_d\cdot |x-y|\cdot(\cM_R |\nabla\varphi|(x)+\cM_R|\nabla \varphi|(y)).\label{Es2}
\end{align}
(iii) Let $\varphi\in L^p_{loc}(\mR^d)$ for some $p>1$. Then for some $C_{d,p}>0$ and any $N,R>0$,
\begin{align}
\left(\int_{B_N}(\cM_R|\varphi|(x))^p\dif x\right)^{1/p}
\leq C_{d,p}\left(\int_{B_{N+R}}|\varphi(x)|^p\dif x\right)^{1/p}.\label{Es30}
\end{align}
\el

\section{SDEs with jumps and smooth coefficients}

In this section, we consider the following SDE with jump:
\begin{align}
X_t(x)&=x+\int^t_0 b_s(X_s(x))\dif s+\int^t_0\sigma_s(X_s(x))\dif W_s
+\int^{t+}_0\!\!\!\int_{\mU}f_s(X_{s-}(x),u)\tilde N(\dif u,\dif s), \label{SDE}
\end{align}
where the coefficients $b:[0,1]\times\mR^d\to\mR^d,\sigma:[0,1]\times\mR^d\to\mR^{d\times d}$ and
$f:[0,1]\times\mR^d\times \mU\to\mR^d$ are measurable functions and smooth in the spatial variable $x$, and satisfy that
\begin{align}
\int^1_0(|b_s(0)|+\|\nabla b_s\|_\infty)\dif s+\int^1_0(|\sigma_s(0)|^2+\|\nabla\sigma_s\|^2_\infty)\dif s<+\infty.\label{Con0}
\end{align}
Moreover, we assume that there exist two functions $L_1, L_2:\mU\to\mR_+$ with
\begin{align}
0\leq L_1(u)\leq \alpha\wedge L_2(u),\ \ \int^1_0\!\!\!\int_{\mU}|L_2(u)|^2(1+L_2(u))^p\nu_s(\dif u)\dif s<+\infty,\label{Con}
\end{align}
where $\alpha\in(0,1)$ is small and $p\in(1,\infty)$ is arbitrary,
and such that for all $(s,x,u)\in[0,1]\times\mR^d\times\mU$,
\begin{align}
|\nabla_x f_s(x,u)|\leq L_1(u), \ \ |f_s(0,u)|\leq L_2(u).\label{CCon}
\end{align}
Under conditions (\ref{Con0})-(\ref{CCon}) with small $\a$ (saying less than $\frac{1}{8d}$),
it is well known that SDE (\ref{SDE}) defines a flow of $C^\infty$-diffeomorphisms
(cf. \cite{Fu-Ku, Pr}, \cite[Theorem 1.3]{Qi-Zh}).

Let $$J_t:=J_t(x):=\nabla X_t(x)\in\mM_{d\times d}.$$
Then $J_t$ satisfies the following SDE (cf. \cite{Fu-Ku, Pr}):
\begin{align}
J_t=\mI+\int^t_0 \nabla b_s(X_s)J_s\dif s+\int^t_0\nabla\sigma_s(X_s)J_s\dif W_s
+\int^{t+}_0\!\!\!\int_{\mU}\nabla f_s(X_{s-},u)J_s\tilde N(\dif u,\dif s).\label{Jac}
\end{align}
The following lemma will be our starting point in the sequel development.
\bl\label{Le5}
The Jacobi determinant $\det(J_t)$ has the following explicit formula:
$$
\det(J_t)=\exp A_t\cdot\exp\left\{M_t-\tfrac{1}{2}\<M^\mathrm{c}\>_t\right\}
\prod_{0<s\leq t}(1+\Delta M_s)e^{-\Delta M_s}=:\exp A_t\cdot \cE(M)_t,
$$
where $A_t:=A^{(1)}_t+A^{(2)}_t$ and $M_t:=M^{\mathrm{c}}_t+M^{\mathrm{d}}_t$ are given by (\ref{Ek1}), (\ref{Ek2}),
(\ref{Ek3}) and (\ref{Ek4}) below.
\el
\begin{proof}
By (\ref{Jac}), It\^o's formula and Lemma \ref{Le1}, we have
\begin{align*}
\det(J_t)&=1+\int^t_0 \div b_s(X_s)\det(J_s)\dif s+\int^t_0\div\sigma_s(X_s)\det(J_s)\dif W_s\no\\
&\quad+\frac{1}{2}\sum_{i,j,k}\int^t_0[\p_i\sigma^{ik}_s\p_j\sigma^{jk}_s-\p_j\sigma^{ik}_s\p_i\sigma^{jk}_s](X_s)
\det(J_s)\dif s\no\\
&\quad+\int^{t+}_0\!\!\!\int_{\mU}\Big[\det((\mI+\nabla f_s(X_{s-},u))J_{s-})-\det(J_{s-})\Big]\tilde N(\dif u,\dif s)\no\\
&\quad+\int^t_0\!\!\!\int_{\mU}\Big[\det((\mI+\nabla f_s(X_{s-},u))J_{s-})-\det(J_{s-})\no\\
&\quad\qquad\qquad-\div f_s(X_{s-},u)\det(J_{s-})\Big]\nu_s(\dif u)\dif s\no\\
&=:1+\int^{t+}_0\det(J_{s-})\dif (A_s+M_s),
\end{align*}
where $A_t:=A^{(1)}_t+A^{(2)}_t$ is a continuous increasing process given by
\begin{align}
A^{(1)}_t=\int^t_0\Big[\div b_s(X_s)+\tfrac{1}{2}\sum_{i,j,k}[\p_i\sigma^{ik}_s\p_j\sigma^{jk}_s
-\p_j\sigma^{ik}_s\p_i\sigma^{jk}_s](X_s)\Big]\dif s\label{Ek1}
\end{align}
and
\begin{align}
A^{(2)}_t=\int^t_0\!\!\!\int_{\mU}\Big[\det(\mI+\nabla f_s(X_{s-},u))
-1-\div f_s(X_{s-},u)\Big]\nu_s(\dif u)\dif s;\label{Ek2}
\end{align}
and $M_t:=M^{\mathrm{c}}_t+M^\mathrm{d}_t$ is a martingale given by
\begin{align}
M^{\mathrm{c}}_t:=\int^t_0\div\sigma_s(X_s) \dif W_s\label{Ek3}
\end{align}
and
\begin{align}
M^\mathrm{d}_t:=\int^{t+}_0\!\!\!\int_{\mU}\Big[\det(\mI+\nabla f_s(X_{s-},u))-1\Big]
\tilde N(\dif u,\dif s).\label{Ek4}
\end{align}
By Dol\'eans-Dade's exponential formula (cf. \cite{Pr}), we obtain the desired formula.
\end{proof}

Below, we shall give an estimate for the $p$-order moment of the Jacobi determinant.
For this aim, we introduce the following function of jump size control $\a$:
\begin{align}
\beta_\a:=(d\a+d!d^2\a^2(1+\a)^{d-2})^{-1}.\label{beta}
\end{align}
Note that
$$
\lim_{\a\downarrow 0}\beta_\a=+\infty.
$$
\bl\label{LL2}
Let $\beta_\a$ be defined by (\ref{beta}), where $\a$ is from (\ref{Con}) small enough so that $\beta_\a>1$. Then
for any $p\in(0,\beta_\a)$, we have
$$
\sup_{x\in\mR^d}\mE\left(\sup_{t\in[0,1]}\det(J_t(x))^{-p}\right)\leq
C\left(p,\int^1_0\|[\div b_s]^-\|_\infty\dif s,\int^1_0\|\nabla\sigma_s\|^2_\infty\dif s,
\int^1_0\!\!\!\int_{\mU}L_1(u)^2\nu_s(\dif u)\dif s\right),
$$
where for a real number $a$, $a^-=\min(-a,0)$,
the constant $C$ is an increasing function with respect to its arguments.
\el
\begin{proof}
First of all, by (\ref{Ek1}), we have
$$
-A^{(1)}_t\leq C\int^1_0(\|[\div b_s]^-\|_\infty+\|\nabla\sigma_s\|^2_\infty)\dif s,
$$
and by (\ref{Ek2}), (\ref{Ep5}) and (\ref{CCon}),
$$
-A^{(2)}_t\leq C\int^t_0\!\!\!\int_{\mU}L_1(u)^2\nu_s(\dif u)\dif s.
$$
Hence, for any $p\geq 0$, we have
$$
\sup_{t\in[0,1]}\exp(-pA_t)\leq
\exp\left(C\int^1_0(\|[\div b_s]^-\|_\infty+\|\nabla\sigma_s\|^2_\infty)\dif s+C\int^1_0\!\!\!\int_{\mU}L_1(u)^2\nu_s(\dif u)\dif s\right).
$$
Thus, by Lemma \ref{Le5}, it suffices to prove that for any $p\in(0,\beta_\a)$,
\begin{align}
\mE\left(\sup_{t\in[0,1]}\cE(M)_t^{-p}\right)\leq C\left(p,
\int^1_0\|\div\sigma_s\|^2_\infty\dif s,\int^1_0\!\!\!\int_{\mU}L_1(u)^2\nu_s(\dif u)\dif s\right).\label{Es1}
\end{align}
Noting that
$$
\Delta M_s:=M_s-M_{s-}=\det(\mI+\nabla f_s(X_{s-},p_s))-1,
$$
by (\ref{Ep5}) and (\ref{Con}), we have
\begin{align}
|\Delta M_s|&\leq|\div f_s(X_{s-},p_s)|+d! d^2L_1(u)^2(1+L_1(u))^{d-2}\no\\
&\leq dL_1(u)+d! d^2L_1(u)^2(1+L_1(u))^{d-2}\label{PP1}\\
&\leq d\a+d!d^2\a^2(1+\a)^{d-2}=\beta_\a^{-1}.\no
\end{align}
Fixing $q\in(p,\beta_\a)$, we also have
$$
|\Delta (-qM)_s|=q|\Delta M_s|<1.
$$
Thus, by Theorem \ref{Th1}, one knows that $t\mapsto\cE(-qM)_t$ is an exponential martingale.
Observe that
\begin{align*}
\cE(M)^{-p}_t&=\cE(-q M)^{\frac{p}{q}}_t\cdot\exp\left\{\frac{(q+1)p}{2}\<M^{\mathrm{c}}\>_t\right\}\cdot\prod_{0<s\leq t}
\frac{(1+\Delta M_s)^{-p}} {(1-q \Delta M_s)^{\frac{p}{q}}}\\
&\leq\cE(-q M)^{\frac{p}{q}}_t\cdot\exp\left\{C\int^1_0\|\div\sigma_s\|^2_\infty\dif s\right\}
\cdot\prod_{0<s\leq t}G(\Delta M_s),
\end{align*}
where
$$
G(r):=\frac{(1+r)^{-p}} {(1-q r)^{\frac{p}{q}}}, \ \ |r|\leq\beta_\a^{-1}.
$$
By H\"older's inequality and Doob's inequality, we obtain that for $\gamma\in(1,\frac{q}{p})$ and $\gamma^*=\frac{\gamma}{\gamma-1}$,
\begin{align}
\mE\left(\sup_{t\in[0,1]}\cE(M)_t^{-p}\right)&\leq C\left(\mE\sup_{t\in[0,1]}\cE(-q M)^{\frac{\gamma p}{q}}_t
\right)^{\frac{1}{\gamma}}
\cdot\left(\mE\prod_{0<s\leq 1}G(\Delta M_s)^{\gamma^*}\right)^{\frac{1}{\gamma^*}}\no\\
&\leq C\left(\mE \cE(-q M)^{\frac{\gamma p}{q}}_1\right)^{\frac{1}{\gamma}}
\cdot\left(\mE\prod_{0<s\leq 1}G(\Delta M_s)^{\gamma^*}\right)^{\frac{1}{\gamma^*}}\no\\
&\leq C\left(\mE\prod_{0<s\leq 1}G(\Delta M_s)^{\gamma^*}\right)^{\frac{1}{\gamma^*}}.\label{Es0}
\end{align}
Thanks to the following limit
$$
\lim_{r\downarrow 0}\frac{\log G(r)}{r^2}=\frac{p(q+1)}{2},
$$
we have for some $C=C(q,p,\beta_\a)>0$,
$$
\left|\log G(r)\right|\leq C|r|^2, \ \ \forall |r|\leq \beta_\a^{-1}.
$$
Therefore, by Lemma \ref{Le3},
\begin{align}
\mE\left[\prod_{0<s\leq 1}G(\Delta M_{s})^{\gamma^*}\right]&=
\mE\exp\left\{\sum_{0<s\leq 1}\gamma^*\log G(\Delta M_{s})\right\}
\leq\mE\exp\left\{\sum_{0<s\leq 1}C|\Delta M_{s}|^2\right\}\no\\
&\stackrel{(\ref{PP1})}{\leq}\mE\exp\left\{\sum_{0<s\leq 1}CL_1(p_s)^2\right\}
\leq\exp\left\{C\int^1_0\!\!\!\int_{\mU}L_1(u)^2\nu_s(\dif u)\dif s\right\}.\label{Es4}
\end{align}
Estimate (\ref{Es1}) now follows by combining (\ref{Es0}) and (\ref{Es4}).
\end{proof}
In order to give an estimate for $\det(\nabla X^{-1}_t(x))$ in terms of $\det(J_t(x))=\det(\nabla X_t(x))$,
we shall use a trick due to Cruzeiro \cite{Cr} (see also \cite{Ci-Cr, Zh2, Fa-Lu-Th}).
Below, let
\begin{align*}
\mu(\dif x):=\frac{\dif x}{(1+|x|^2)^d}.
\end{align*}
We write
$$
\cJ_t(\omega,x):=\frac{(X_t(\omega,\cdot))_\sharp\mu(\dif x)}{\mu(\dif x)},\ \
\cJ^-_t(\omega,x):=\frac{(X^{-1}_t(\omega,\cdot))_\sharp\mu(\dif x)}{\mu(\dif x)},
$$
which means that for any nonnegative measurable function $\varphi$ on $\mR^d$,
\begin{align}
\int_{\mR^d}\varphi(X_t(\omega,x))\mu(\dif x)&=\int_{\mR^d}\varphi(x) \cJ_t(\omega,x)\mu(\dif x),\label{P0}\\
\int_{\mR^d}\varphi(X^{-1}_t(\omega, x))\mu(\dif x)&=\int_{\mR^d}\varphi(x) \cJ^-_t(\omega, x)\mu(\dif x).\label{P1}
\end{align}
It is easy to see that for almost all $\omega$ and all $(t,x)\in[0,1]\times\mR^d$,
\begin{align}
\cJ_t(\omega,x)=[\cJ^-_t(\omega,X^{-1}_t(\omega,x))]^{-1}\label{P2}
\end{align}
and
\begin{align}
\cJ^-_t(x)=\frac{(1+|x|^2)^d}{(1+|X_t(x)|^2)^d}\det(J_t(x)).\label{P3}
\end{align}

We need the following estimate:
\bl\label{LL1}
For any $p\geq 1$, we have
\begin{align*}
&\sup_{x\in\mR^d}\mE\left(\sup_{t\in[0,1]}\frac{(1+|X_t(x)|^2)^p}{(1+|x|^2)^p}\right)\\
&\qquad\leq C\left(p,
\int^1_0\left\|\frac{|b_s(x)|}{1+|x|}\right\|_{\infty}\dif s,
\int^1_0\left\|\frac{|\sigma_s(x)|}{1+|x|}\right\|_{\infty}^2\dif s,\int^1_0\!\!\!\int_{\mU}L_2(u)^2(1+L_2(u))^{4p-2}\nu_s(\dif u)\dif s\right),
\end{align*}
where the constant $C$ is an increasing function with respect to its arguments.
\el
\begin{proof}
Letting $h(x):=(1+|x|^2)^p$, by It\^o's formula, we have
\begin{align*}
h(X_t)-h(x)&=\int^t_0(b^i_s\p_ih)(X_s)\dif s
+\int^t_0(\sigma^{ik}_s\p_ih)(X_s)\dif W^k_s+\frac{1}{2}\int^t_0(\p_i\p_jh\cdot\sigma^{ik}_s\sigma^{jk}_s)(X_s)\dif s\\
&\quad+\int^t_0\!\!\!\int_{\mU}(h(X_{s-}+f_s(X_{s-},u))-h(X_{s-})
-f^i_s(X_{s-},u)\p_ih(X_{s-}))\nu_s(\dif u)\dif s\\
&\quad+\int^{t+}_0\!\!\!\int_{\mU}(h(X_{s-}+f_s(X_{s-},u))-h(X_{s-})) \tilde N(\dif u,\dif s).
\end{align*}
By elementary calculations, one has
\begin{align}
C_1(1+|x|)^{2p}\leq h(x)\leq C_2(1+|x|)^{2p}\label{PP3}
\end{align}
and
\begin{align*}
|\p_ih(x)|\leq \frac{Ch(x)}{1+|x|}\leq C(1+|x|)^{2p-1},
\ \ |\p_i\p_jh(x)|\leq \frac{Ch(x)}{(1+|x|)^2}\leq C(1+|x|)^{2p-2}.
\end{align*}
On the other hand, by Taylor's formula, we have
$$
|h(x+y)-h(x)|\leq |y^i\p_ih(x+\theta_1 y)|
$$
and
$$
|h(x+y)-h(x)-y^i\p_ih(x)|\leq |y^iy^j\p_i\p_jh(x+\theta_2 y)|/2,
$$
where $\theta_1,\theta_2\in(0,1)$.
Thus, for $p\geq 1$, we have
\begin{align*}
|h(x+f_s(x,u))-h(x)|&\leq |f_s(x,u)|\cdot(1+|x+\theta_1 f_s(x,u)|)^{2p-1}\\
&\stackrel{(\ref{CCon})}{\leq} (L_2(u)+L_1(u)|x|)(1+L_2(u)+(1+L_1(u))|x|)^{2p-1}\\
&\stackrel{(\ref{Con})}{\leq} L_2(u)(1+L_2(u))^{2p-1}(1+|x|)^{2p}\\
&\stackrel{(\ref{PP3})}{\leq} L_2(u)(1+L_2(u))^{2p-1}h(x)
\end{align*}
and
$$
|h(x+f_s(x,u))-h(x)-f^i_s(x,u)\p_ih(x)|\leq L_2(u)^2(1+L_2(u))^{2p-2}h(x).
$$
Using the above estimates, if we let
$$
\ell_1(s):=\left\|\frac{|b_s(x)|}{1+|x|}\right\|_{\infty},
\ \ \ell_2(s):=\left\|\frac{|\sigma_s(x)|}{1+|x|}\right\|_{\infty},
\ \ \ell_3(s):=\int_\mU L_2(u)^2(1+L_2(u))^{4p-2}\nu_s(\dif u),
$$
then, by Burkholder's inequality and Young's inequality, we have
\begin{align*}
\mE\left(\sup_{s\in[0,t]}h(X_s)\right)&\leq h(x)+C\int^t_0(\ell_1(s)+\ell_2^2(s))\mE h(X_s)\dif s
+\mE\left(\int^t_0\ell_2^2(s)h(X_s)^2\dif s\right)^{1/2}\\
&\quad +C\int^t_0\ell_3(s) \mE h(X_s)\dif s+ C\mE\left(\int^t_0\ell_3(s)h(X_s)^2\dif s\right)^{1/2}\\
&\leq h(x)+C\int^t_0(\ell_1(s)+\ell_2^2(s)+\ell_3(s))\mE h(X_s)\dif s+\frac{1}{2}\mE\left(\sup_{s\in[0,t]}h(X_s)\right),
\end{align*}
which leads to
$$
\mE\left(\sup_{s\in[0,t]}h(X_s)\right)\leq h(x)+C\int^t_0(\ell_1(s)+\ell_2^2(s)+\ell_3(s))\mE h(X_s)\dif s.
$$
Hence, by Gronwall's inequality, we obtain
$$
\mE\left(\sup_{s\in[0,1]}h(X_s)\right)\leq C h(x).
$$
The proof is complete.
\end{proof}

Combining Lemmas \ref{LL2} and \ref{LL1}, we obtain that
\bt\label{Th0}
Let $\beta_\a$ be defined by (\ref{beta}), where $\a$ is from (\ref{Con}) small enough so that $\beta_\a>1$. Then
for any $p\in (0,\beta_\a)$,
$$
\mE\left(\sup_{t\in[0,1]}\int_{\mR^d}|\cJ_t(x)|^{p+1}\mu(\dif x)\right)\leq C,
$$
where the constant $C$ is inherited from Lemmas \ref{LL2} and \ref{LL1}.
\et
\begin{proof}
The estimate follows from
$$
\int_{\mR^d}|\cJ_t(x)|^{p+1}\mu(\dif x)\stackrel{(\ref{P2})(\ref{P1})}{=}\int_{\mR^d}|\cJ^-_t(x)|^{-p}\mu(\dif x)
\stackrel{(\ref{P3})}{=}\int_{\mR^d}\frac{(1+|X_t(x)|^2)^{dp}}{(1+|x|^2)^{dp}}\det(J_t(x))^{-p}\mu(\dif x),
$$
H\"older's inequality and Lemmas \ref{LL2} and \ref{LL1}.
\end{proof}

\section{SDEs with jumps and rough drifts}

We first introduce the following notion of generalized stochastic flows (cf. \cite{Le-Li1,Zh1,Zh2}).
\bd\label{Def1}
Let $X_t(\omega,x)$ be a $\mR^d$-valued measurable
stochastic field on $[0,1]\times\Omega\times\mR^d$. For a locally finite measure $\mu$ on $\mR^d$,
we say $X$ a {\bf $\mu$-almost everywhere stochastic flow} or {\bf generalized stochastic flow} of SDE (\ref{SDE})  if
\begin{enumerate}[{\bf(A)}]
\item for some $p\geq 1$, there exists a constant $K_p>0$
such that for any nonnegative measurable function $\varphi\in L^p_\mu(\mR^d)$,
\begin{align}\label{Den}
\sup_{t\in[0,1]}\mE\int_{\mR^d}\varphi(X_t(x))\mu(\dif x)\leq K_p\|\varphi\|_{L^p_\mu};
\end{align}
\item for  $\mu$-almost all $x\in\mR^d$, $t\mapsto X_t(x)$ is a c\'adl\'ag and
($\sF_t$)-adapted  process and solves equation (\ref{SDE}).
\end{enumerate}
\ed

The main result of this section is:
\bt\label{Main}
Assume that for some $q>1$,
$$
|\nabla b|\in L^1([0,1]; L^q_{loc}(\mR^d)),\ \ [\div b]^-, |\nabla\sigma|^2,
\frac{|b|}{1+|x|},\frac{|\sigma|^2}{1+|x|^2} \in L^1([0,1]; L^\infty(\mR^d)),
$$
and for some functions $L_i:\mU\to[0,+\infty), i=1,2$ satisfying (\ref{Con}), and all $(s,u)\in[0,1]\times\mU$ and $x,y\in\mR^d$,
\begin{align}
|f_s(x,u)-f_s(y,u)|\leq L_1(u)|x-y|, \ \ |f_s(0,u)|\leq L_2(u).\label{SI}
\end{align}
Let $\mu(\dif x)=(1+|x|^2)^{-d}\dif x$ and let $\beta_\a$ be defined by (\ref{beta}),
where $\a$ is from (\ref{Con}) small enough so that $\beta_\a>\frac{1}{q-1}$. Then there exists a unique $\mu$-almost
everywhere stochastic flow to SDE (\ref{SDE}) with any $p\geq q$ in (\ref{Den}).
\et
\br
Let $b(x)=\frac{x}{|x|}1_{x\not=0}$. It is easy to check that $\div b(x)=\frac{d-1}{|x|}$ and $|\nabla b|\in L^{p}_{loc}(\mR^d)$ provided
that $p\in[1,d)$.
\er

Let $\chi\in C^\infty(\mR^d)$ be a nonnegative cutoff function with
\begin{equation} \label{eq:4}
\|\chi\|_\infty\leq 1,\ \
\chi(x)=\left\{ \begin{aligned}
&1,\ \ |x|\leq 1, \\
&0,\ \ |x|\geq 2.
\end{aligned} \right.
\end{equation}
Let $\rho\in C^\infty(\mR^d)$ be a nonnegative mollifier with support in $B_1:=\{|x|\leq 1\}$ and $\int_{\mR^d}\rho(x)\dif x=1$.
Set
$$
\chi_n(x):=\chi(x/n),\ \ \rho_n(x):=n^{d}\rho(nx)
$$
and define
\begin{align}
b_s^n:=b_s*\rho_n\cdot\chi_n,\ \
\sigma^n_s:=\sigma_s*\rho_n,\ \ f^n_s(\cdot,u)=f_s(\cdot,u)*\rho_n.\label{Es3}
\end{align}
The following lemma is direct from the definitions and the property of convolutions.
\bl\label{Le9}
For some $C>0$ independent of $n$, we have
$$
\int^1_0\|[\div b_s^n]^-\|_\infty\dif s\leq\int^1_0\|[\div b_s]^-\|_\infty\dif s
+C\int^1_0\left\|\frac{|b_s(x)|}{1+|x|}\right\|_\infty\dif s
$$
$$
\int^1_0\left\|\frac{|b^n_s(x)|}{1+|x|}\right\|_\infty\dif s
\leq C\int^1_0\left\|\frac{|b_s(x)|}{1+|x|}\right\|_\infty\dif s
$$
$$
\int^1_0\|\nabla\sigma^n_s\|_\infty^2\dif s\leq \int^1_0\|\nabla\sigma_s\|_\infty^2\dif s,\\
$$
$$
\int^1_0\left\|\frac{|\sigma^n_s(x)|}{1+|x|}\right\|_\infty^2\dif s
\leq C\int^1_0\left\|\frac{|\sigma_s(x)|}{1+|x|}\right\|^2_\infty\dif s
$$
and
$$
|\nabla_x f^n_s(x,u)|\leq L_1(u),\ \ |f^n_s(0,u)|\leq 2L_2(u).
$$
\el
\begin{proof}
The first estimate follows from that
$$
\div b^n_s(x)=(\div b_s)*\rho_n(x)\cdot\chi_n(x)+b^i_s*\rho_n(x)\cdot\p_i\chi_n(x)
$$
and
$$
|\p_i\chi_n(x)|=\frac{|(\p_i\chi)(x/n)|}{n}\leq\frac{C1_{n\leq |x|\leq n+1}}{1+|x|}.
$$
The other estimates are similar.
\end{proof}

Let $X^n_t(x)$ be the stochastic flow of $C^\infty$-diffeomorphisms to SDE (\ref{SDE}) associated with coefficients
$(b^n,\sigma^n,f^n)$.

\bl\label{Le0}
let $\beta_\a$ be defined by (\ref{beta}),
where $\a$ is from (\ref{Con}) small enough so that $\beta_\a>1$. Then for any $p>1+\frac{1}{\beta_\a}$,
there exists a constant $C_p>0$ such that for all non-negative function $\varphi\in L^p_\mu(\mR^d)$,
\begin{align}
\sup_{n}\mE\left(\sup_{t\in[0,1]}\int_{\mR^d}\varphi(X^n_t(x))\mu(\dif x)\right)\leq C_p\|\varphi\|_{L^p_\mu}.\label{BB4}
\end{align}
\el
\begin{proof}
The estimate follows from
\begin{align*}
\int_{\mR^d}&\varphi(X^n_t(x))\mu(\dif x)=
\int_{\mR^d}\varphi(x) \cJ^n_t(x)\mu(\dif x) \leq\|\varphi\|_{L^p_\mu}\left(\int_{\mR^d}
|\cJ^n_t(x)|^{\frac{p}{p-1}}\mu(\dif x)\right)^{1-\frac{1}{p}},
\end{align*}
and Theorem \ref{Th0} and Lemma \ref{Le9}.
\end{proof}

\bl\label{Le-1}
For any $n,m>4/\delta>0$, we have
$$
\frac{|z+f^n_t(x+z,u)-f^m_t(x,u)|^2-|z|^2}{|z|^2+\delta^2}\leq 4(L_1(u)+L_1(u)^2).
$$
\el
\begin{proof}
Noticing that by the property of convolutions and (\ref{SI}),
\begin{align*}
|f^n_t(x+z,u)-f^m_t(x,u)|&\leq |f^n_t(x+z,u)-f^n_t(x,u)|+|f^n_t(x,u)-f_t(x,u)|
+|f^m_t(x,u)-f_t(x,u)|\\
&\leq L_1(u)|z|+L_1(u)(n^{-1}+m^{-1}),
\end{align*}
we have
\begin{align*}
\frac{|z+f^n_t(x+z,u)-f^m_t(x,u)|^2-|z|^2}{|z|^2+\delta^2}&\leq\frac{2|z|(|z|+(n^{-1}+m^{-1}))L(u)
+(|z|+(n^{-1}+m^{-1}))^2L_1(u)^2}{|z|^2+\delta^2}\\
&\leq 2[1+\delta^{-1}(n^{-1}+m^{-1})](L_1(u)+L_1(u)^2),
\end{align*}
which yields the desired estimate.
\end{proof}

We now prove the following key estimate.
\bl\label{Le8}
For any $R>1$, there exist constants $C_1,C_2>0$ such that for all $\delta\in(0,1)$ and $n,m>4/\delta$,
\begin{align}
&\mE\int_{G^{n,m}_R}
\sup_{t\in[0,1]}\log\left(\frac{|X^n_t(x)-X^m_t(x)|^2}{\delta^2}+1\right)\mu(\dif x)\no\\
&\qquad\leq C_1 +\frac{C_2}{\delta}\int^1_0\left(\|b^n_s-b^m_s\|_{L^q(B_{R})}
+\|\sigma^n_s-\sigma^m_s\|^2_{L^{2q}(B_{R})}\right)\dif s,\label{Es8}
\end{align}
where $\mu(\dif x)=(1+|x|^2)^{-d}\dif x$ and
$G^{n,m}_R(\omega):=\Big\{x\in\mR^d: \sup_{t\in[0,1]}|X^n_t(\omega,x)|\vee|X^m_t(\omega,x)|\leq R\Big\}$.
\el
\begin{proof}
Set
\begin{align*}
Z^{n,m}_t(\omega,x):=X^n_t(\omega,x)-X^m_t(\omega, x)
\end{align*}
and
$$
F^{n,m}_t(\omega,x,u):=f^n_t(X^n_{t-}(\omega,x),u)- f^m_t(X^m_{t-}(\omega,x),u).
$$
If there are no confusions, we shall drop the variable ``$x$'' below. Note that
\begin{align*}
Z^{n,m}_t=\int^t_0(b^n_s(X^n_s)-b^m_s(X^m_s))\dif s+\int^t_0(\sigma^n_s(X^n_s)-\sigma^m_s(X^m_s))\dif W_s
+\int^{t+}_0\!\!\!\int_{\mU}F^{n,m}_s(u)\tilde N(\dif u,\dif s).
\end{align*}
By It\^o's formula, we have
\begin{align*}
\log\left(\frac{|Z^{n,m}_t|^2}{\delta^2}+1\right)
&=2\int^t_0\frac{\<Z^{n,m}_s, b^n_s(X^n_s)-b^m_s(X^m_s)\>}{|Z^{n,m}_s|^2+\delta^2}\dif s
+2\int^t_0\frac{\<Z^{n,m}_s, (\sigma^n_s(X^n_s)-\sigma^m_s(X^m_s))\dif W_s\>}{|Z^{n,m}_s|^2+\delta^2}\\
&\quad+\int^t_0\frac{\|\sigma^n_s(X^n_s)-\sigma^m_s(X^m_s)\|^2}{|Z^{n,m}_s|^2+\delta^2}\dif s
-2\int^t_0\frac{|(\sigma^n_s(X^n_s)-\sigma^m_s(X^m_s))^\mathrm{t}\cdot Z^{n,m}_s|^2}{(|Z^{n,m}_s|^2+\delta^2)^2}\dif s\\
&\quad+\int^{t+}_0\!\!\!\int_{\mU}\left(\log\frac{|Z^{n,m}_{s-}+F^{n,m}_s(u)|^2+\delta^2}
{|Z^{n,m}_{s-}|^2+\delta^2}-\frac{|Z^{n,m}_{s-}+F^{n,m}_s(u)|^2-|Z^{n,m}_{s-}|^2}
{|Z^{n,m}_{s-}|^2+\delta^2}\right)\nu_s(\dif u)\dif s\\
&\quad+\int^{t+}_0\!\!\!\int_{\mU}\log\frac{|Z^{n,m}_{s-}+F^{n,m}_s(u)|^2+\delta^2}
{|Z^{n,m}_{s-}|^2+\delta^2}\tilde N(\dif u,\dif s)\\
&=:I^{n,m}_1(t)+I^{n,m}_2(t)+I^{n,m}_3(t)+I^{n,m}_4(t)+I^{n,m}_5(t)+I^{n,m}_6(t).
\end{align*}
For $I^{n,m}_1(t)$,  we have
\begin{align*}
\sup_{t\in[0,1]}|I^{n,m}_1(t)|&\leq2\int^1_0\frac{|b^n_s(X^n_s)-b^n_s(X^m_s)|}{\sqrt{|Z^{n,m}_s|^2+\delta^2}}\dif s
+\frac{2}{\delta}\int^1_0|b^n_s(X^m_s)-b^m_s(X^m_s)|\dif s=:I^{n,m}_{11}+I^{n,m}_{12}.
\end{align*}
Noting that
$$
G^{n,m}_R(\omega)\subset \{x: |X^n_t(\omega,x)|\leq R\}\cap\{x: |X^m_t(\omega,x)|\leq R\},\ \ \forall t\in[0,1],
$$
we have
\begin{align}
\mE\int_{G^{n,m}_R}|I^{n,m}_{12}(x)|\mu(\dif x)&\leq
\frac{2}{\delta}\mE\int^1_0\!\!\!\int_{\mR^d}|1_{B_R}(b^n_s-b^m_s)|(X^m_s(x))\mu(\dif x)\dif s\no\\
&\stackrel{(\ref{BB4})}{\leq}\frac{C}{\delta}\int^1_0\|1_{B_R}(b^n_s-b^m_s)\|_{L^q_\mu}\dif s\no\\
&\leq\frac{C}{\delta}\int^1_0\|b^n_s-b^m_s\|_{L^q(B_R)}\dif s.\label{PL0}
\end{align}
For $I^{n,m}_{11}$, in view of $\mu(\dif x)\leq\dif x$, we have
\begin{align}
\mE\int_{G^{n,m}_R}|I^{n,m}_{11}(x)|\mu(\dif x)&\stackrel{(\ref{Es2})}{\leq} C\mE\int^1_0\!\!\!\int_{G^{n,m}_R}(\cM_{2R}|\nabla b^n_s|(X^n_s(x))
+\cM_{2R}|\nabla b^n_s|(X^m_s(x)))\mu(\dif x)\dif s\no\\
&\stackrel{(\ref{BB4})}{\leq} C\int^1_0\left(\int_{B_R}(\cM_{2R}|\nabla b^n_s|(x))^q\mu(\dif x)\right)^{1/q}\dif s\no\\
&\stackrel{(\ref{Es30})}{\leq} C \int^1_0\|\nabla b^n_s\|_{L^q(B_{3R})}\dif s\leq C \int^1_0\|\nabla b_s\|_{L^q(B_{3R})}\dif s.\label{PL1}
\end{align}
For $I^{n,m}_2(t)$, set
$$
\tau^{n,m}_R(\omega,x):=\inf\Big\{t\in[0,1]: |X^n_t(\omega,x)|\vee X^m_t(\omega,x)>R\Big\},
$$
then
$$
G^{n,m}_R(\omega)=\{x:\tau^{n,m}_R(\omega,x)=1\}.
$$
By Burkholder's inequality and Fubini's theorem, we have
\begin{align*}
&\mE\int_{G^{n,m}_R}\sup_{t\in[0,1]}|I^{n,m}_{2}(t,x)|\mu(\dif x)\\
&\leq\int_{\mR^d}\mE\left(\sup_{t\in[0,\tau^{n,m}_R(x)]}\left|\int^{t}_0\frac{\<Z^{n,m}_s(x),(\sigma^n_s(X^n_s(x))
-\sigma^m_s(X^m_s(x)))\dif W_s\>}{|Z^{n,m}_s(x)|^2+\delta^2}\right|\right)\mu(\dif x)\\
&\leq C\int_{\mR^d}\mE\left[\int^{ \tau^{n,m}_R(x)}_0\frac{|Z^{n,m}_s(x)|^2|\sigma^n_s(X^n_s(x))-\sigma^m_s(X^m_s(x))|^2}
{(|Z^{n,m}_s(x)|^2+\delta^2)^2}\dif s\right]^{\frac{1}{2}}\mu(\dif x)\\
&\leq C\mu(\mR^d)^{\frac{1}{2}}\left[\mE\int^1_0\!\!\!\int_{G^{n,m}_R}
\frac{|\sigma^n_s(X^n_s(x))-\sigma^m_s(X^m_s(x))|^2}
{|Z^{n,m}_s(x)|^2+\delta^2}\mu(\dif x)\dif s\right]^{\frac{1}{2}}.
\end{align*}
As the treatment of $I^{n,m}_1(t)$, by Lemma \ref{Le0}, we can prove that
\begin{align}
\mE\int_{G^{n,m}_R}\sup_{t\in[0,1]}|I^{n,m}_{2}(t,x)|\mu(\dif x)
\leq \left(C\int^1_0\|\nabla\sigma_s\|^2_{L^{2q}(B_{R+1})}\dif s
+\frac{C}{\delta}\int^1_0\|\sigma^n_s-\sigma^m_s\|^2_{L^{2q}(B_{R})}\dif s\right)^{\frac{1}{2}},\label{PL2}
\end{align}
and similarly,
\begin{align}
\mE\int_{G^{n,m}_R}\sup_{t\in[0,1]}|I^{n,m}_{3}(t,x)|\mu(\dif x)
&\leq C\int^1_0\|\nabla\sigma_s\|^2_{L^{2q}(B_{R+1})}\dif s
+\frac{C}{\delta}\int^1_0\|\sigma^n_s-\sigma^m_s\|^2_{L^{2q}(B_{R})}\dif s.\label{PL33}
\end{align}
Since $I^{n,m}_4(t)$ is negative, we can drop it. For $I^{n,m}_5(t)$, by Lemma \ref{Le-1} and the elementary inequality
$$
|\log(1+r)-r|\leq C|r|^2,\ \ r\geq -\tfrac{1}{2},
$$
we have
\begin{align}
\mE\int_{G^{n,m}_R}\sup_{t\in[0,1]}|I^{n,m}_5(t,x)|\mu(\dif x)
&\leq C\int^1_0\!\!\!\int_{\mU}L_1(u)^2\nu_s(\dif u)\dif s.\label{PL3}
\end{align}
For $I^{n,m}_6(t)$, as in the treatment of $I^{n,m}_2(t)$ and $I^{n,m}_5(t)$, we also have
\begin{align}
\mE\int_{G^{n,m}_R}\sup_{t\in[0,1]}|I^{n,m}_6(t,x)|\mu(\dif x)
&\leq C\left(\int^1_0\!\!\!\int_{\mU}L_1(u)^2\nu_s(\dif u)\dif s\right)^{\frac{1}{2}}.\label{PL4}
\end{align}
Combining (\ref{PL0})-(\ref{PL4}), we obtain (\ref{Es8}).
\end{proof}

We are now in a position to give

\begin{proof}[Proof of Theorem \ref{Main}]
Set
$$
\Phi^{n,m}(x):=\sup_{t\in[0,1]}|X^n_t(x)-X^m_t(x)|
$$
and
$$
\Psi^{n,m}_\delta(x):=\log\left(\frac{\Phi^{n,m}(x)^{2}}{\delta^2}+1\right).
$$
We have
\begin{align*}
\mE\int_{\mR^d}\Phi^{n,m}(x)\mu(\dif x)
&=\mE\int_{(G^{n,m}_R)^c}\Phi^{n,m}(x)\mu(\dif x)
+\mE\int_{G^{n,m}_R}\Phi^{n,m}(x)\mu(\dif x),
\end{align*}
where $G^{n,m}_R$ is defined as in Lemma \ref{Le8}.
By Lemmas \ref{LL1} and \ref{Le9}, the first term is less than
$$
\frac{1}{\sqrt{R}}\mE\int_{\mR^d}\left(\sup_{t\in[0,1]}|X^n_t(x)|^{\frac{3}{2}}
+\sup_{t\in[0,1]}|X^m_t(x)|^{\frac{3}{2}}\right)\mu(\dif x)\leq
\frac{C\int_{\mR^d}(1+|x|^2)^{\frac{3}{4}-d}\dif x}{\sqrt{R}}\leq\frac{C}{\sqrt{R}},
$$
where $C$ is independent of $n,m$ and $R$, and $d\geq 2$.

For the second term, we make the following decomposition:
\begin{align*}
\mE\!\!\int_{G^{n,m}_R}\!\!\!\Phi^{n,m}(x)\mu(\dif x)=
\mE\!\!\int_{G^{n,m}_R\cap\{\Psi^{n,m}_\delta\geq\eta\}}\!\!\!\Phi^{n,m}(x)\mu(\dif x)
+\mE\!\!\int_{G^{n,m}_R\cap\{\Psi^{n,m}_\delta<\eta\}}\!\!\!\Phi^{n,m}(x)\mu(\dif x)
=:I^{n,m}_1+I^{n,m}_2.
\end{align*}
For $I^{n,m}_1$, by H\"older's inequality, Lemma \ref{LL1} and (\ref{Es8}), we have
\begin{align*}
I^{n,m}_1\leq\frac{C_R}{\sqrt{\eta}}\left(\mE\int_{G^{n,m}_R}\!\!\!\Psi^{n,m}_\delta(x)\mu(\dif x)\right)^{\frac{1}{2}}
\leq \frac{C_R}{\sqrt{\eta}}+\frac{C_R}{\sqrt{\delta\eta}}\left(\int^1_0\!\!\!\left(\|b^n_s-b^m_s\|_{L^q(B_{R})}
+\|\sigma^n_s-\sigma^m_s\|^2_{L^{2q}(B_{R})}\right)\dif s\right)^{\frac{1}{2}}.
\end{align*}
For $I^{n,m}_2$, noticing that if $\Psi^{n,m}_\delta(x)\leq\eta$, then $\Phi^{n,m}(x)\leq\delta\sqrt{e^\eta-1}$, we have
$$
I^{n,m}_2\leq C\delta\sqrt{e^\eta-1}.
$$
Combining the above calculations, we obtain that
$$
\mE\int_{\mR^d}\!\!\!\Phi^{n,m}(x)\mu(\dif x)
\leq \frac{C}{\sqrt{R}}+\frac{C_R}{\sqrt{\eta}}+C\delta e^{\eta/2}+
\frac{C_R}{\sqrt{\delta\eta}}\left(\int^1_0\!\!\!\left(\|b^n_s-b^m_s\|_{L^q(B_{R})}
+\|\sigma^n_s-\sigma^m_s\|^2_{L^{2q}(B_{R})}\right)\dif s\right)^{\frac{1}{2}}.
$$
Taking limits in order: $n,m\to\infty$, $\delta\to 0$, $\eta\to\infty$  and $R\to\infty$ yields that
\begin{align*}
\lim_{n,m\to\infty}\mE\int_{\mR^d}\left(\sup_{t\in[0,1]}|X^n_t(x)-X^m_t(x)|\right)\mu(\dif x)=0.
\end{align*}
Thus, there exists an adapted c\'adl\'ag process $X_t(x)$ such that
$$
\lim_{n\to\infty}\mE\int_{\mR^d}\left(\sup_{t\in[0,1]}|X^n_t(x)-X_t(x)|\right)\mu(\dif x)=0.
$$
By Lemma \ref{Le4}, it is standard to check that $X_t(x)$ solves SDE (\ref{SDE}) in the sense of Definition \ref{Def1}.

For the uniqueness, let $X^i_t(x), i=1,2$ be two almost everywhere stochastic flows of SDE (\ref{SDE}). As in the proof of Lemma \ref{Le8}, we have
\begin{align*}
\mE\int_{G_R}\sup_{t\in[0,1]}\log\left(\frac{|X^1_t(x)-X^2_t(x)|^2}{\delta^2}+1\right)\mu(\dif x)\leq C
\end{align*}
where $G_R(\omega):=\Big\{x\in\mR^d: \sup_{t\in[0,1]}|X^1_t(\omega,x)|\vee|X^2_t(\omega,x)|\leq R\Big\}$ and $C$ is independent of $\delta$.
Letting $\delta\to 0$ and $R\to\infty$, we obtain that $X^1_t(\omega, x)=X^2_t(\omega, x)$ for all $t\in[0,1]$ and $P\times\mu$-almost all $(\omega,x)$.
\end{proof}

\section{Probabilistic representation for the solutions of PIDEs}

In this section we work in the canonical space $\Omega=\mD^d_{[0,1]}$:
the set of all right continuous functions with left limits. The generic element in $\Omega$ is denoted by $w$.
The space $\Omega$ can be endowed with two complete metrics:
uniform metric and Skorohod metric. We remark that only under Skorohod metric, $\Omega$
is separable. For $t\in[0,1]$, let $\sF_t:=\sigma\{w_s: s\in[0,t]\}$ and
set $\sF=\sF_1$. Then $\sF$ coincides with the $\sigma$-algebra generated by Skorohod's topology. For a Polish space $E$,
by $\cP(E)$ we denote the space of all Borel probability measures over $E$.

Below we consider the more general L\'evy generator:
$$
\sL_t\varphi(x):=\frac{1}{2} a^{ij}_t(x)\p_i\p_j \varphi(x)+b^i_t(x)\p_i \varphi(x)+
\int_{\mR^d\setminus\{0\}}\left[\varphi(x+y)-\varphi(x)-\frac{\<y,\nabla\varphi(x)\>}{1+|y|^2}\right]\nu_t(\dif y),
$$
where $a^{ij}_t(x):=\sum_k\sigma^{ik}_t(x)\sigma^{jk}_t(x)$ and $\nu$ satisfies that
$$
\int^1_0\!\!\!\int_{\mR^d\setminus\{0\}}\frac{|y|^2}{1+|y|^2}\nu_t(\dif y)\dif t<+\infty.
$$

We recall the following notion of Stroock and Varadhan's martingale solutions (cf. \cite{St, St-Va}).
\bd
(Martingale Solutions)
Let $\mu_0\in\cP(\mR^d)$. A probability measure $P$ on $(\Omega,\sF)$ is called a martingale solution
corresponding to the operator $\sL$ and initial law $\mu_0$ if $\mu_0=P\circ w_0^{-1}$ and
for all $\varphi\in C^\infty_0(\mR^d)$,
$$
\varphi(w_t)-\varphi(w_0)-\int^t_0(\sL_s\varphi)(w_s)\dif s
$$
is a $P$-martingale with respect to $(\sF_t)$, which is equivalent that for all $\theta\in\mR^d$,
\begin{align*}
&\exp\Bigg[\mathrm{i}\<\theta, w_t-w_0-\int^t_0b_s(w_s)\dif s\>-\frac{1}{2}\int^t_0a^{ij}_s(w_s)\theta^i\theta^j\dif s\\
&\qquad-\int^t_0\!\!\!\int_{\mR^d\setminus\{0\}}\left(e^{\mathrm{i}\<\theta,y\>}
-1-\frac{\mathrm{i}\<\theta,y\>}{1+|y|^2}\right)\nu_s(\dif y)\dif s\Bigg]
\end{align*}
is a $P$-martingale with respect to $(\sF_t)$.
\ed
For any $w\in\Omega$ and $\Gamma\in\cB(\mR^d\setminus\{0\})$, we define
$$
\eta(t,w,\Gamma):=\sum_{0<s\leq t}1_{\Gamma}(w(s)-w(s-))
$$
and
$$
\tilde\eta(t,w,\Gamma):=\eta(t,w,\Gamma)-\int^t_0\nu_s(\Gamma)\dif s.
$$
The following result is from \cite[Corollaries 1.3.1 and 1.3.2]{St}.
\bt\label{Th2}
Let $P\in\cP(\Omega)$ be a martingale solution corresponding to $(\sL,\mu_0)$. Given $\delta>0$, define
$$
\gamma^\delta_t(w):=w_t-\int_{|y|<\delta}y\tilde\eta(t,w,\dif y)-\int_{|y|\geq \delta}y\eta(t,w,\dif y)
$$
and
\begin{align}
\hat b^\delta_t(x):=b_t(x)+\int_{|y|<\delta}\frac{y|y|^2}{1+|y|^2}\nu_t(\dif y)-\int_{|y|\geq\delta}\frac{y}{1+|y|^2}\nu_t(\dif y).\label{BB1}
\end{align}
Then $M(t,w):=\gamma^\delta_t(w)-\int^t_0\hat b^\delta_s(w_s)\dif s$ is independent of $\delta>0$ and continuous and $(\sF_t)$-adapted.
Moreover, for any $\theta\in\mR^d$ and $|g(y)|^2\leq\frac{C|y|^2}{1+|y|^2}$,
$$
t\mapsto \exp\left[\mathrm{i}\<\theta,M(t)-M(0)\>+\frac{1}{2}\int^t_0a^{ij}_s(w_s)\theta^i\theta^j\dif s\right]
$$
and
$$
t\mapsto\int_{\mR^d\setminus\{0\}}g(y)\tilde\eta(t,w,\dif y)
$$
are $P$-martingales with respect to $(\sF_t)$.
\et
Let us now consider the following integro-differential equation of Fokker-Planck type:
\begin{align}
\p_t\mu_t=\sL^*_t\mu_t,\label{PDE1}
\end{align}
where $\sL^*_t$ is the formal adjoint operator of $\sL_t$ given by
$$
\sL^*_t\mu:=\frac{1}{2} \p_i\p_j (a^{ij}_t(x)\mu)-\p_i(b^i_t(x)\mu)+
\int_{\mR^d\setminus\{0\}}\left[\tau_{y}\mu-\mu+\frac{y^i\p_i\mu}{1+|y|^2}\right]\nu_t(\dif y).
$$
Here, PIDE (\ref{PDE1}) is understood in the distributional sense, i.e., for any $\varphi\in C^\infty_b(\mR^d)$,
\begin{align}
\p_t\<\mu_t, \varphi\>=\<\mu_t,\sL_t \varphi\>.\label{BB2}
\end{align}
If $\mu_t(\dif x)=u_t(x)\dif x$, then (\ref{PDE1}) reads as
\begin{align}
\p_t u_t=\sL^*_t u_t.\label{PDE11}
\end{align}

The following result gives the uniqueness of measure-valued solutions for (\ref{PDE1}) in the case of smooth coefficients.
\bt\label{Th4}
Assume that $a$ and $b$ are smooth and satisfies that for all $k\in\{0\}\cup\mN$,
$$
\sup_{t\in[0,1]}\|\nabla^k a^{ij}_t\|_\infty+\sup_{t\in[0,1]}\|\nabla^k b^i_t\|_\infty<+\infty.
$$
Then for any $\mu_0\in\cP(\mR^d)$, PIDE (\ref{PDE1})
admits a unique measure-valued solution $\mu_t\in\cP(\mR^d)$.
\et
\begin{proof}
The existence is clear as introduced in the introduction. Let us now prove the uniqueness.
 For $0\leq s<t\leq 1$ and $x\in\mR^d$,
let $X_{s,t}(x)$ solve the following SDE:
$$
X_{s,t}(x)=x+\int^t_s\hat b_r(X_{s,r}(x))\dif r+\int^t_s\sqrt{a_r}(X_{s,r}(x))\dif W_r
+\int_{B_1^0}y\tilde N(\dif y,(s,t])+\int_{B^c_1}yN(\dif y,(s,t]),
$$
where $\hat b_r(x)$ is defined by (\ref{BB1}) with $\delta=1$, $\sqrt{a_r}$ denotes the square root of symmetric nonnegative matrix $a_r$
and $N(\dif y,\dif t)$ is a Poisson random point measure with intensity measure $\nu_t(\dif y)\dif t$, $B^0_1:=B_1\setminus\{0\}$ and
$B^c_1=\mR^d\setminus B_1$.
For any $\varphi\in C^\infty_b(\mR^d)$, define
$$
\cT_{s,t}\varphi(x):=\mE(\varphi(X_{s,t}(x))).
$$
Then $\cT_{s,t}\varphi(x)\in C^\infty_b(\mR^d)$ and for all $0\leq s<r<t\leq 1$,
$$
\cT_{s,r}\cT_{r,t}\varphi(x)=\cT_{s,t}\varphi(x).
$$
It is easy to verify that
$$
\p_s\cT_{s,t}\varphi+\sL_s\cT_{s,t}\varphi=0.
$$
Let $\mu^i_t, i=1,2$ be two solutions of PIDE (\ref{PDE1}) with the same initial values.
Then by (\ref{BB2}), we have
$$
\p_s\<\mu^i_s,\cT_{s,t}\varphi\>=\<\mu^i_s,\p_s\cT_{s,t}\varphi+\sL_s\cT_{s,t}\varphi\>=0,\ \ i=1,2.
$$
Since $\mu^1_0=\mu^2_0$, we have
$$
\<\mu^1_s,\cT_{s,t}\varphi\>=\<\mu^2_s,\cT_{s,t}\varphi\>, \ \ s\in[0,t].
$$
In particular,
$$
\<\mu^1_t,\varphi\>=\<\mu^2_t,\varphi\>,
$$
which implies that $\mu^1_t=\mu^2_t$ for any $t\in[0,1]$.
\end{proof}
We now prove the following extension of Figalli's result \cite[p.116, Theorem 2.6]{Fi}, which is originally due to
Ambrosio \cite{Am}.
\bt\label{Th3}
Assume that $b$ and $a$ are bounded and measurable functions.
Let $\mu_t\in\cP(\mR^d)$ be a measure-valued solution of PIDE (\ref{PDE1}) with initial value $\mu_0\in\cP(\mR^d)$. Then
there exists a martingale solution $P\in\cP(\Omega)$ corresponding to $(\sL,\mu_0)$
such that for all $t\in[0,1]$ and $\varphi\in C^\infty_0(\mR^d)$,
\begin{align}
\<\mu_t,\varphi\>=\mE^P(\varphi(w_t)).\label{PLO1}
\end{align}
\et
\begin{proof}
Let $\rho: \mR^d \to(0,+\infty)$ be a convolution kernel such that $|\nabla^k\rho (x)|\leq C_k\rho(x)$
for any $k\in\mN$ (for instance $\rho(x)=e^{-|x|^2/2}/(2\pi)^{d/2}$). Let $\rho_\eps(x):=\eps^{-d}\rho(x/\eps)$, $\eps>0$,
and define
$$
\mu^\eps_t:=\mu_t*\rho_\eps,\ \ b^\eps_t:=\frac{(b_t\mu_t)*\rho_\eps}{\mu^\eps_t},
\ \ a^\eps_t:=\frac{(a_t\mu_t)*\rho_\eps}{\mu^\eps_t}.
$$
It is easy to see that for any $k\in\{0\}\cup\mN$,
$$
\|\nabla^k b^\eps_t\|_\infty\leq C_k\|\nabla^k b_t\|_\infty,\ \
\|\nabla^k a^\eps_t\|_\infty\leq C_k\|\nabla^k a_t\|_\infty.
$$
With a little abuse of notation, we are denoting the measure $\mu^\eps_t$ and its density with respect to the Lebesgue measure by the same symbol.
If we take the convolutions with $\rho_\eps$ for both sides of PIDE (\ref{PDE1}), then
$$
\p_t\mu^\eps_t=\frac{1}{2}\p_i\p_j(a^{\eps,ij}_t\mu^\eps_t)-\p_i(b^{\eps,i}_t\mu^\eps_t)+\int_{\mR^d\setminus\{0\}}
\left[\tau_y\mu^\eps_t-\mu^\eps_t+\frac{\<y,\nabla\mu^\eps_t\>}{1+|y|^2}\right]\nu_t(\dif y),
$$
subject to $\mu^\eps_0=\mu_0*\rho_\eps$.
By Theorem \ref{Th4}, the unique solution to this PIDE can be represented by
$$
\mu^\eps_t=\mbox{Law of }X^\eps_t,
$$
i.e., for any $\varphi\in C^\infty_0(\mR^d)$,
\begin{align}
\<\mu^\eps_t,\varphi\>=\mE\varphi(X^\eps_t),\label{PLK}
\end{align}
where $X^\eps_t$ solves the following SDE with jump
$$
X^\eps_t=X^\eps_0+\int^t_0\hat b^\eps_s(X^\eps_s)\dif s+\int^t_0\sqrt{a^{\eps}_s}(X^\eps_s)\dif W_s
+\int_{B^0_1}y\tilde N(\dif y,(0,t])+\int_{B^c_1}yN(\dif y,(0,t]),
$$
and the law of $X^\eps_0$ is $\mu^\eps_0$. Here, $\hat b^\eps_s(x)$ is defined by (\ref{BB1}) with $\delta=1$ and replacing $b$ by $b^\eps$.

Let $P_\eps$ be the law of $t\mapsto X^\eps_t$ in $\Omega$. Since
$$
P_\eps(|w_0|\geq R)=\mu_0^\eps(B_R^c)\to 0\mbox{ uniformly in $\eps$ as $R\to\infty$},
$$
by \cite[p.237, Theorem A.1]{St}, $(P_\eps)_{\eps\in(0,1)}$ is tight in $\cP(\Omega)$. Let $P$ be any accumulation of $(P_\eps)_{\eps\in(0,1)}$.
Without loss of generality, we assume that $P_\eps$ weakly converges to $P$ as $\eps\to 0$. By taking weak limits
for both sides of (\ref{PLK}), it is clear that (\ref{PLO1}) holds.

For completing the proof, it remains to show that $P$ is a martingale solution corresponding to $(\sL,\mu_0)$. That is, we need to prove that
for any $0\leq s<t\leq 1$ and bounded continuous and $\sF_s$-measurable function $\Phi^s$ on $\Omega$, $\varphi\in C^\infty_0(\mR^d)$,
$$
\mE^P\left[\left(\varphi(w_t)-\varphi(w_s)-\int^t_s(\sL_r\varphi)(w_r)\dif r\right)\Phi^s(w)\right]=0.
$$
This will follow by taking weak limits for
$$
\mE^{P_\eps}\left[\left(\varphi(w_t)-\varphi(w_s)-\int^t_s(\sL^\eps_r\varphi)(w_r)\dif r\right)\Phi^s(w)\right]=0.
$$
The more details can be found in \cite[p.118, Step 3]{Fi}.
\end{proof}
\bd
(Weak solution) If there exists a filtered probability space $(\Omega,\sF,P;(\sF_t)_{t\in[0,1]})$
and an $(\sF_t)$-adapted Brownian motion $W_t$, an $(\sF_t)$-adapted Poisson random measure $N(\dif y,\dif t)$ with
intensity measure $\nu_t(\dif y)\dif t$
and an $(\sF_t)$-adapted process $X_t$ on $(\Omega,\sF,P;(\sF_t)_{t\in[0,1]})$ such that for some $\delta>0$ and all $t\in[0,1]$,
\begin{align}
X_t&=X_0+\int^t_0\hat b^\delta_s(X_s)\dif s+\int^t_0\sigma_s(X_s)\dif W_s
+\int_{B^0_\delta}y\tilde N(\dif y,(0,t])+\int_{B^c_\delta}yN(\dif y,(0,t]),
\end{align}
where $\hat b^\delta_s(x)$ is defined by (\ref{BB1}),
then we say $(\Omega,\sF,P;(\sF_t)_{t\in[0,1]})$ together with $(W,N,X)$ a weak solution.
By weak uniqueness, we means that any two weak solutions
with the same initial law have the same law in $\Omega$.
\ed

The following result gives the equivalence between weak solutions and martingale solutions.

\bt\label{Th5}
The existence of martingale solutions implies the existence of weak solutions. In particular,
the uniqueness of weak solutions implies the uniqueness of martingale solutions.
\et
\begin{proof}
Let $P\in\cP(\Omega)$ be a martingale solution.
By Theorem \ref{Th2}, one knows that under $P$, $\eta$ is a Poisson random point measure
with intensity measure $\nu_t(\dif y)\dif t$ and $M$ is a continuous martingale with covariation process
$$
\<M^i,M^j\>_t=\frac{1}{2}\sum_k\int^t_0(\sigma^{ik}_s\sigma^{jk}_s)(w_s)\dif s.
$$
Let $(\hat\Omega,\hat\sF,\hat P; (\hat\sF_t)_{t\in[0,1]})$ be another filtered probability space supporting a Brownian motion
$\hat W_t$. Let $(\tilde\Omega,\tilde\sF,\tilde P; (\tilde\sF_t)_{t\in[0,1]})$ be the product filtered probability space of
$(\Omega,\sF,P; (\sF_t)_{t\in[0,1]})$ and $(\hat\Omega,\hat\sF,\hat P; (\hat\sF_t)_{t\in[0,1]})$.
Let $\pi:\tilde\Omega\to\Omega$ be the canonical projection. Define
$$
\tilde M_t(\tilde\omega):=M_t(\pi(\tilde\omega)),\ \ \tilde \sigma_t(\tilde\omega):=\sigma_t(\pi(\tilde\omega)_t)
$$
and
$$
\tilde N_t(\tilde\omega,\dif y):=\eta(t,\pi(\tilde\omega),\dif y),\ \ \tilde X_t(\tilde\omega):=\pi(\tilde\omega).
$$
Then by the proof of \cite[p.108, Theorem 4.5.2]{St}, there exists another Brownian motion $(\tilde W_t)_{t\in[0,1]}$ defined on
$(\tilde\Omega,\tilde\sF,\tilde P; (\tilde\sF_t)_{t\in[0,1]})$ such that
$$
\tilde M_t=\int^t_0\tilde \sigma_s\dif \tilde W_s,\ \ \tilde P-a.s.
$$
Hence, $(\tilde\Omega,\tilde\sF,\tilde P; (\tilde\sF_t)_{t\in[0,1]})$ together with
$(\tilde W,\tilde N,\tilde X)$ is a weak solution.
\end{proof}

The main result of this section is:
\bt
Assume that for some $q>1$,
$$
|\nabla b|\in L^\infty([0,1]; L^q_{loc}(\mR^d;\mR^d)),\ \ [\div b]^-,
|b|,|\sigma|, |\nabla\sigma|\in L^\infty([0,1]\times\mR^d),
$$
and for any $p\geq 1$,
$$
\int^1_0\!\!\!\int_{\mR^d\setminus\{0\}}|y|^2(1+|y|)^p\nu_t(\dif y)\dif t<+\infty.
$$
Let $r>\frac{q}{q-1}=q^*$.
Then for any probability density function $\phi$ with
$$
\int_{\mR^d}\phi(x)^r(1+|x|^2)^{(r-1)d}\dif x<+\infty,
$$
there exists a unique solution $u_t$ to PIDE (\ref{PDE11}) in the class of
$$
\sM_{q^*}:=\left\{u_t\in L^{q^*}(\mR^d): u_t(x)\geq 0, \int_{\mR^d} u_t(x)\dif x=1,
\sup_{t\in[0,1]}\int_{\mR^d}u_t(x)^{q^*}(1+|x|^2)^{(q^*-1)d}\dif x<+\infty\right\}.
$$
Moreover, if $q>d$, then the uniqueness holds in the measure-valued space $\cP(\mR^d)$.
\et
\begin{proof}
(Existence) Set $\mu(\dif x):=\dif x/(1+|x|^2)^d$ and
let $X_t(x)$ be the $\mu$-almost everywhere stochastic flow of the following SDE
$$
X_t(x)=x+\int^t_0 \hat b_s(X_s(x))\dif s+\int^t_0\sigma_s(X_s(x))\dif W_s
+\int_{\mR^d\setminus\{0\}}y\tilde N(\dif y,(0,t]),
$$
where $N(\dif y,(0,t])$ is a Poisson random measure with intensity $\nu_t(\dif y)\dif t$ and
$$
\hat b_s(x):=b_s(x)+\int_{\mR^d\setminus\{0\}}\frac{y|y|^2}{1+|y|^2}\nu_s(\dif y).
$$
Since in this case, $L_1=0$ in Theorem \ref{Main}, the $p$ in (\ref{Den}) can be arbitrarily close to $1$.
Let $X_0$ be an $\sF_0$-measurable random variable with law $\phi(x)\dif x$.  Define
$$
Y_t:=X_t(X_0).
$$
It is easy to check that $Y_t$ solves the following SDE:
\begin{align}
Y_t=X_0+\int^t_0 \hat b_s(Y_s)\dif s+\int^t_0\sigma_s(Y_s)\dif W_s+\int_{\mR^d\setminus\{0\}}y\tilde N(\dif y,(0,t]).\label{BB5}
\end{align}
Now for any $\varphi\in C^\infty_0(\mR^d)$, by H\"older's inequality, we have
\begin{align*}
\mE\varphi(Y_t)&=\mE(\mE\varphi(X_t(x))|x=X_0)=\int_{\mR^d}\mE\varphi(X_t(x)) \phi(x)\dif x\\
&\leq\left(\int_{\mR^d}|\mE\varphi(X_t(x))|^{\frac{r}{r-1}}\mu(\dif x)\right)^{1-\frac{1}{r}}\left(\int_{\mR^d}
\Big(\phi(x)(1+|x|^2)^d\Big)^r\mu(\dif x)\right)^{\frac{1}{r}}\\
&=\left(\mE\int_{\mR^d}|\varphi(X_t(x))|^{\frac{r}{r-1}}\mu(\dif x)\right)^{1-\frac{1}{r}}\left(\int_{\mR^d}
\phi(x)^r(1+|x|^2)^{(r-1)d}\dif x\right)^{\frac{1}{r}}\stackrel{(\ref{Den})}{\leq} C_\phi\|\varphi\|_{L^q_\mu},
\end{align*}
which then implies that $Y_t$ has an absolutely continuous probability density $u_t\in\sM_{q^*}$ with
$$
\int_{\mR^d}u_t(x)\varphi(x)\dif x=\mE\varphi(Y_t)\leq C_\phi\|\varphi\|_{L^q_\mu},\ \ \frac{1}{q^*}+\frac{1}{q}=1.
$$
By It\^o's formula, it is immediate that $u_t$ solves PIDE (\ref{PDE11}) in the distributional sense.

(Uniqueness) Let $u_t^i\in\sM_{q^*}$ be any two solutions of PIDE (\ref{PDE11}) with the same initial value $u_0=\phi$.
Let $P^i\in\cP(\Omega)$ be two  martingale solutions corresponding to $\mu^i_t(\dif x)=u^i_t(x)\dif x$ by Theorem \ref{Th3}.
Since  for any $\varphi\in C_0^\infty(\mR^d)$,
$$
\int_{\mR^d}u^i_t(x)\varphi(x)\dif x=\mE^{P^i}\varphi(w_t),\ \ i=1,2,
$$
we only need to prove that $P^1=P^2$. By Theorem \ref{Th5} and \cite[p.104, Theorem 137]{Si}, it suffices to prove the pathwise
uniqueness of SDE (\ref{BB5}). Let $Y^i_t, i=1,2$ be two solutions of SDE (\ref{BB5}) defined on the same filtered probability space supporting
a Brownian motion $W$ and a Poisson random measure $N$ with intensity measure $\nu_t(\dif y)\dif t$,
where for any $\varphi\in C_0^\infty(\mR^d)$,
$$
\int_{\mR^d}u^i_t(x)\varphi(x)\dif x=\mE\varphi(Y^i_t),\ \ \ i=1,2.
$$
Since $u^i_t\in\sM_p$, by suitable approximation, we have for any $\varphi\in L^q_\mu(\mR^d)$,
\begin{align}
\sup_{t\in[0,1]}\mE\varphi(Y^i_t)\leq C\|\varphi\|_{L^q_\mu},\ \ i=1,2.\label{BB6}
\end{align}
Set
$$
Z_t:=Y^1_t-Y^2_t,\ \ \ \tau_R:=\inf\{t\in[0,1]: |Y^1_t|\vee|Y^2_t|>R\}.
$$
Basing on (\ref{BB6}), as in the proof of Lemma \ref{Le8},  we have for any $\delta>0$,
\begin{align}
\mE\log\left(\frac{|Z_{t\wedge\tau_R}|^2}{\delta^2}+1\right)
&\leq 2\mE\int^{t\wedge\tau_R}_0\frac{\<Z_s, b_s(Y^1_s)-b_s(Y^2_s)\>}{|Z_s|^2+\delta^2}\dif s
+\mE\int^{t\wedge\tau_R}_0\frac{\|\sigma_s(Y^1_s)-\sigma_s(Y^2_s)\|^2}{|Y^1_s-Y^2_s|^2+\delta^2}\dif s\no\\
&\stackrel{(\ref{Es2})}{\leq} C\mE\int^{t\wedge\tau_R}_0\Big(\cM_{2R}|\nabla b_s|(Y^1_s)+\cM_{2R}|\nabla b_s|(Y^2_s)\Big)\dif s
+\int^1_0\|\nabla\sigma_s\|_\infty^2\dif s\label{E1}\\
&\leq C\int^t_0\|1_{B_R}\cdot \cM_{2R}|\nabla b_s|\|_{L^q_\mu}\dif s+\int^1_0\|\nabla\sigma_s\|_\infty^2\dif s\no\\
&\leq C\int^t_0\|\cM_{2R}|\nabla b_s|\|_{L^q(B_R)}\dif s+\int^1_0\|\nabla\sigma_s\|_\infty^2\dif s\no\\
&\stackrel{(\ref{Es30})}{\leq}C\int^1_0\|\nabla b_s\|_{L^q(B_{3R})}\dif s
+\int^1_0\|\nabla\sigma_s\|_\infty^2\dif s,\no
\end{align}
where $C$ is independent of $\delta$.
Letting first $\delta\to 0$ and then $R\to\infty$, we obtain that $Z_t=0$ a.s., i.e., $Y^1_t=Y^2_t$ a.s.

In the case of $q>d$, let $Y^1_t$ be the solution constructed in the proof of existence and $Y^2_t$ another solution of SDE (\ref{BB5})
corresponding to any measure-valued solution $\mu_t$ with $\mu_0(\dif x)=\phi(x)\dif x$. In the above proof of (\ref{E1}), 
instead of using (\ref{Es2}), we use Morrey's inequality
(\ref{Es02}) to deduce that $Y^1_t=Y^2_t$.
\end{proof}

{\bf Acknowledgements:}

The supports of NSFs of China (No. 10971076; 10871215) are acknowledged.

\end{document}